\documentclass[english, 12pt]{article}

\usepackage{amsfonts, amsmath}
\usepackage{amssymb}
\usepackage{amscd}
\usepackage{amscd}
\usepackage{amsthm}

\newtheorem{theo}{Theorem}[section]
\newtheorem{prop}{Proposition}[section]                   
\newtheorem{lem}{Lemma}[section]

\numberwithin{equation}{section}

\def\R{\mathbb{R}}
\def\N{\mathbb{N}}

\newcommand{\rhup}{\rightharpoonup}
\newcommand{\rc}{C_{\rm Riesz}}
\newcommand{\soc}{C_{\rm sob}}
\newcommand{\Lb}{L^2}
\newcommand{\Lbr}{L^{\rho'}}
\newcommand{\W}{W^{1,2}}
\newcommand{\Wd}{W^{-1,2}}

\newcommand{\varp}{\varphi}

\newcommand{\alp}{\alpha}

\begin{document}
\title{Local existence and uniqueness for the frictional Newton-Schr\"odinger
  equation in three dimensions}

\author{\normalsize
Ali BenAmor\\{\em I.P.E.I. Tunis}\\
{\em 2, Rue Jawaher Lel Nehru}\\
{\em 1008, Tunis}\\
{\em Tunisia}\\
ali.benamor@ipeit.rnu.tn
\and  Philippe Blanchard\\{\em Fakult\"at f{\"u}r Physik}\\
{\em Uni. Bielefeld}\\
{\em D-33501 Bielefeld}\\
{\em Germany}\\
blanchard@physik.uni-bielefeld.de}

\date{}
\maketitle

\begin{abstract} We  prove, in this paper, local existence and
  uniqueness of solution for the frictional Newton-Schr\"odinger
 equation in three dimensions. Further we show that the blow-up alternative
 holds true as well as the continuous dependence of the solution w.r.t. the initial
 data. Our method is rather direct  and based
 essentially on a fixed point-type theorem due to Weissinger and an
 approximation process.

\end{abstract}


\section{Introduction}

We consider the frictional Newton-Schr\"odinger equation (frNSE for
short) in three dimensions:
\begin{eqnarray}
\begin{cases}
i\frac{\partial\psi}{\partial t}=-\frac{\hbar}{2m}\Delta\psi+i4\pi\frac{Gm^2}{\hbar}\psi\int_{\R^3}\frac{|\psi(y)|^2}{|\cdot-y|}\,dy-i4\pi\frac{Gm^2}{\hbar}\psi\int_{\R^3}\int_{\R^3}\frac{|\psi(y)|^2|\psi(z)|^2}{|z-y|}\,dy\,dz & \\
\psi(0)=\varphi &  \end{cases},
\label{frnse}
\end{eqnarray} 

where $\varphi$ is a given element from the Sobolev space
$W^{1,2}:=W^{1,2}(\R^3)$.\\
This equation is the limit case as $R\to 0$ of the alternative (frNSE)
considered by Diosi \cite{diosi07}. Indeed he considered a (frNSE) with kernel $U(x,x')$ in the second and third
term of (\ref{frnse}) that behaves like

\begin{eqnarray}
U(x,x')\sim
\begin{cases}
c+\frac{1}{2}m\omega_G^2|x-x'|^2,\  {\rm if}\ |x-x'|<<R  & \\
-Gm^2|x-x'|^{-1},\  {\rm if}\ |x-x'|>>R &  \end{cases},
\end{eqnarray}
where $R>0$.\\
As observed by Diosi, equation (\ref{frnse}) differs from the usual
Newton-Schr\"odinger equation (NSE for short)
\cite{diosi07,diosi84,penrose,tod01,tod99,adler}:
\begin{eqnarray}
\begin{cases}
i\frac{\partial\psi}{\partial t}=-\frac{\hbar}{2m}\Delta\psi-4\pi\frac{Gm^2}{\hbar}\psi\int_{\R^3}\frac{|\psi(y)|^2}{|\cdot-y|}\,dy& \\
\psi(0)=\varphi &  \end{cases},
\label{nse}
\end{eqnarray}
essentially by the term $-i$, which is responsible for the friction, and also the
last term (for normalization).\\

By the way, we would like to mention that the (NSE) in $1$, $2$ and $3$
dimensions has a long standing history. It was already investigated in the
50's by Pekar \cite[pp.29-34]{pekar}. Since that time, there is a huge literature dealing
with the equation, we cite, as examples, the papers of Lieb \cite{lieb},
Penrose \cite{penrose1}, Nawa-Ozawa \cite{ozawa-nonlocal}, the
very instructive book of Cazenave \cite{cazenave} (and references therein),
\cite{kato89} and recently \cite{miao07}.\\
The problem of existence of bound states for the $1$-dimensional (NSE) was
already investigated in a recent paper of Choquard and Stubbe
\cite{choquard}.\\
Also the pseudo-relativistic (NSE) was the subject of papers by Fr\"ohlich et al. \cite{F1,F2}.\\ 
However, to our best knowledge, the frictional Newton-Schr\"odinger equation
was not treated in the literature.\\ 
Unlike the (NSE), the (frNSE) can not have a
stationary solution and has no energy. Further, due the occurrence of complex
coefficients, one can not expect the conservation of charge for equation (\ref{frnse}). Therefore all proof-methods based on
'energy functional' and 'conservation laws'  arguments do not work any more to
prove existence of solutions in this situation.\\
Also, since the nonlinearity is nonlocal, Kato's method \cite{kato87,kato2} (based essentially on
the use of a fixed point theorem and on Strichartz estimates) does not work; but
maybe its generalization given in \cite[p.98]{cazenave}.\\
Here we shall  propose an alternative method which, for proving local
existence and continuous dependence does not rely on Strichartz estimates, and
hence applicable also for the equation with space variable lying in a  subset
of $\R^3$.\\ 
Our main goal, in this paper, is to prove local existence and uniqueness of (weak)-solution for the (frNSE) within
the space $X:=L^{\infty}\big(I,W^{1,2}\big)$, where $I$ is an interval
of the real line containing $0$. We will also show that the blow-up
alternative holds true as well as the continuous dependence w.r.t. the initial
data of the solution.\\
To reach our purposes we will use the following strategy: Truncate the Newton
kernel by eliminating the singularity lying on the diagonal. This consists to
introduce the sequence of integral operators
$$
K_n\phi:=\int_{\{|\cdot-y|>1/n\}}|\cdot-y|^{-1}\phi(y)\,dy.
$$
We obtain in this manner a new function representing the interaction given by
(up to a complex factor) $f_n(\phi)=\phi K_n(|\phi|^2)$ that preserves the space
$\W$. In this stage, using Duhamel's formula together with a fixed point
argument  (Satz of Weissinger) we prove existence and uniqueness of solution
for the approximate problem (replace $\phi K(|\phi|^2)$ by $f_n(\phi)$). The main
ingredient of the proof, in this step, is the local uniform Lipschitz property
of the $f_n$'s.\\  
After this step we prove that the approximate solutions tend to the solution
of  (\ref{frnse}). This will rely upon the crucial property that functions
$g_1,g_2$ (defined below) are weakly continuous.\\
Finally continuous dependence w.r.t. initial data will be established.

\section{Preparatory results}
  
The aim of this section is the proof of some preparatory results for the local
existence and uniqueness theorem.\\
We set $W^{-1,2}$ the dual space of $\W$ and introduce the functions 
\begin{eqnarray}
g_1:\ \W\to\Wd,\
\psi\mapsto\,\psi\int_{\R^3}\frac{|\psi(y)|^2}{|\cdot-y|}\,dy,
\label{g1}
\end{eqnarray}    
and
\begin{eqnarray}
g_2:\ \W\to\Wd,\
\psi\mapsto\,\psi\int_{\R^3}\int_{\R^3}\frac{|\psi(x)|^2|\psi(y)|^2}{|\cdot-y|}\,dy\,dx,
\label{g2}
\end{eqnarray}    
Later we shall prove that $g_1,g_2$  are well-defined, continuous and are bounded on bounded sets.\\
For a fixed $\varphi\in\W$, by a weak solution of (\ref{frnse}) we mean a function 
$$
\psi\in L^{\infty}\big(I,W^{1,2}\big)\cap
W^{1,\infty}\big(I,W^{-1,2}\big),
$$
that satisfies satisfies 
\begin{eqnarray}
i\frac{\partial\psi}{\partial
  t}=-\frac{\hbar}{2m}\Delta\psi+i4\pi\frac{Gm^2}{\hbar}g_1(\psi)-i4\pi\frac{Gm^2}{\hbar}g_2(\psi)\,
  {\rm in}\ W^{-1,2}\, {\rm for\ a.e.}\, t\in I,\, \psi(0)=\varphi.
\label{weak}
\end{eqnarray}

For the convenience of the reader we reproduce the comments made by
Cazenave (see \cite[pp.56-57]{cazenave}) to visualize that the position of the
problem, in this manner is coherent.\\
From the already indicated properties of both functions $g_1,g_2$, we observe 
that if function $\psi\in X:=L^{\infty}\big(I,W^{1,2}\big)$ then both
$g_1(\psi),g_2(\psi)$ are in the space $L^{\infty}\big(I,\Wd\big)$, yielding
$$
-\frac{\hbar}{2m}\Delta\psi+i4\pi\frac{Gm^2}{\hbar}g_1(\psi)-i4\pi\frac{Gm^2}{\hbar}g_2(\psi)\in
L^{\infty}\big(I,\Wd\big).
$$
Thus if $\psi$ satisfies the first part of (\ref{weak}) then it is in
$W^{1,\infty}\big(I,\Wd\big)$.\\
Furthermore, the fact that $\psi\in X\cap
W^{1,\infty}\big(I,\Wd\big)$ yields that $\psi\in C(I,L^2(\R^3))$, which
implies that the second identity in (\ref{weak}) is meaningful.\\
From now on we set
\begin{eqnarray*}
\alp_1:=\frac{\bar h}{2m},\  \alp_2:=4\pi G\frac{m^2}{\bar
  h}.
\end{eqnarray*}
and for every $\psi\in\W$
\begin{eqnarray*}
K\psi:=\int_{\R^3}\frac{\psi(y)}{|\cdot-y|}\,dy,\ G_1(\psi):=\int_{\R^3}|\psi(x)|^2K(|\psi(x)|^2)\,dx\
\end{eqnarray*} 

Let us recall some classical inequalities that we shall use many
times. First, the Sobolev inequality: For every $p\in [2,6]$ there is
a constant $C_{\rm sob}(p)$ such that
\begin{eqnarray}   
\big(\int_{\R^3}|\psi|^pd\,x\big)^{2/p}\leq C_{\rm
  sob}(p)\big(\int_{\R^3}|\nabla\psi|^2d\,x+\int_{\R^3}|\psi|^2d\,x\big),\
  \forall\,\psi\in\W.
\label{sobolev}
\end{eqnarray}
Second the inequality satisfied by Riesz potentials (see \cite{saloff}):\\
For every $1<p<3/2$, set $q:=\frac{3p}{3-2p}$. Then there is a constant
$C_{\rm Riesz}(p)$ such that
\begin{eqnarray}
\|K\psi\|_{L^q}\leq C_{\rm Riesz}(p)\|\psi\|_{L^p},\ \forall\,\psi\in L^p.
\label{riesz}
\end{eqnarray}

As a notation we shall designate by $\int \cdots\,dx$ the integral on $\R^3$
w.r.t. Lebesgue measure  of a given function.


\begin{lem} The functions $g_1,g_2$ are well defined. Moreover they are
  bounded on bounded sets.
\label{defined}
\end{lem}

\begin{proof}  For every $\varp,\psi\in\W$. Making use of inequality
  (\ref{sobolev}) together with (\ref{riesz}) we get
\begin{eqnarray} 
\int|\varp g_1(\psi)|\,dx&\leq& \big(\int|\varp\psi|^{\frac{5}{4}}\big)^{\frac{4}{5}}\big(\int(K(|\psi|^2))^5\big)^{\frac{1}{5}}\leq
 \big(\int|\varp|^{\frac{5}{2}}\big)^{\frac{2}{5}}\big(\int|\psi|^{\frac{5}{2}}\big)^{\frac{2}{5}}\cdot\rc
 \big(\int|\psi|^{\frac{30}{13}}\big)^{\frac{13}{30}}\nonumber\\
&&\leq\rc\soc\soc'\|\varp\|_{\W}\|\psi\|_{\W}^2.
\end{eqnarray}    
Hence for every fixed  $\psi\in\W$, the function
\begin{eqnarray*}
T:=\W\to\R,\ \varp\mapsto {\rm Re}(\int\varp g_1(\psi)\,dx) 
\end{eqnarray*}
 is linear and continuous. Thus for every $\psi\in\W$, $g_1(\psi)\in
 \Wd$ and\\ $\|T=g_1(\psi)\|_{\Wd}\leq\rc\soc\soc'\|\psi\|_{\W}^3$ yielding that $g_1$ is
 bounded on bounded sets.\\
For $g_2$, we have $g_2(\psi)=\psi G_1(\psi)$ and by the same inequalities
\begin{eqnarray} 
G_1(\psi)&&\leq\big(\int|\psi|^{\frac{5}{2}}\big)^{\frac{4}{5}}\cdot\rc
 \big(\int|\psi|^{\frac{30}{13}}\big)^{\frac{13}{30}}\nonumber\\
&&\leq\rc\soc\|\psi\|_{\W}^3.
\end{eqnarray}   
Thus $g_2$ is well defined, bounded on bounded sets and is even in $\W$. 
\end{proof}

The functions $g_1,g_2$ enjoy further interesting properties especially the
local Lipschitz property:

\begin{lem} The function $g_1$ satisfies the following properties: There is $r_1,\rho\in
  [2,6)$ such that
\begin{itemize}
\item[i)] $g_1\in C(\W,L^{\rho'})$.
\item[ii)]  For every $0<M<\infty$, there is a constant $C(M)$ such
  that
$$
\|g_1(\varp)-g_1(\psi)\|_{L^{\rho'}}\leq C(M)\|\varp-\psi\|_{L^{r_1}},
$$  
for every $\varp,\psi\in\W$ such that $\|\varp\|_{\W}\leq M,\
\|\psi\|_{\W}\leq M$.
\end{itemize} 
\label{lip-g1}
\end{lem}

\begin{proof} We have for every $\psi\in\W$, $g_1(\psi)=\psi
  K(|\psi|^2)$ and $\psi\in L^p$ for every $2\leq p\leq 6$.\\
i):  We must first show that there is $\rho\in[2,6)$ such that
\begin{eqnarray*}
\forall\,\psi\in\W,\ g(\psi)\in L^{\rho'}.
\end{eqnarray*}
Let $\rho\in [2,6)$ be fixed and $r$ a real number such that
\begin{eqnarray}
(1):\,1<r<\frac{3}{3-\rho'},\ (2):\,2\leq r\rho'\leq 6\, {\rm and}\,
   (3):\,1<r<\frac{12}{6+\rho'}   
\label{exponent} 
\end{eqnarray} 

Note that if $\psi\in L^{r\rho'}$ and $ K(|\psi|^2)\in L^{r'\rho'}$,
then by H\"older inequality $\psi K(|\psi|^2)\in L^{\rho'}$. Therefore
we will prove that if  $\psi\in\W$ and if $r,\rho$ satisfy the above
conditions then (i) is fulfilled as well as (ii) with $r_1=r\rho'$.\\
Now if  $\psi\in\W$ and if $r,\rho$ satisfy (1-2) of (\ref{exponent}) then
$r'\rho'>3$. Thus setting  $p=\frac{3r'\rho'}{3+2r'\rho'}$ then
$1<p<3/2$. Activating inequality (\ref{riesz}) together with inequality (\ref{sobolev}) gives
\begin{eqnarray*}
\big(\int(K(|\psi|^2))^{r'\rho'}\,dx\big)^{\frac{1}{r'\rho'}}\leq
  \rc\big(\int|\psi|^{2p}\,dx\big)^{\frac{1}{p}}\leq \rc\soc^2\|\psi\|_{\W}^2.
\end{eqnarray*} 

Next from H\"older inequality, we achieve
\begin{eqnarray} 
\big(\int|g_1(\psi)|^{\rho'}\,dx\big)^{\frac{1}{\rho'}}&\leq&\big(\int|\psi|^{r\rho'}\big)^{\frac{1}{r\rho'}}\big(\int(K(|\psi|^2))^{r'\rho'}\big)^{\frac{1}{r'\rho'}}\nonumber\\
&&\leq\soc'\rc\soc^2\|\psi\|_{\W}^3,
\end{eqnarray}

yielding that $g_1:\W\to L^{\rho'}$ is well defined.\\
{\em Continuity}: Let $(\psi_k)\subset\W$ converging in $\W$ to
$\psi$. Then
\begin{eqnarray*}
\int|g_1(\psi_k)-g_1(\psi)|^{\rho'}\,dx\leq 2^{\rho'-1}\int|(\psi_k-\psi)K(|\psi_k|^2)|^{\rho'}\,dx
+2^{\rho'-1}\int|\psi K(|\psi_k|^2-|\psi|^2)|^{\rho'}\,dx.
\end{eqnarray*}
The preceeding calculus shows that 
\begin{eqnarray}
\int|g_1(\psi_k)-g_1(\psi)|^{\rho'}\,dx&\leq&
C\|\psi_k-\psi\|_{\W}\|\psi_k\|_{\W}^2+C\|\psi\|_{\W}\||\psi_k|^2-|\psi|^2\|_{L^p}\nonumber\\
&&\to 0\ {\rm as}\ k\to\infty,
\end{eqnarray}
and $g_1$ is continuous.\\

ii) Fix $M>0$ and $\varp,\psi\in\W$ such that $\|\varp\|_{\W},\|\psi\|_{\W}\leq
M$;  Let $r,\rho$ be as in (\ref{exponent}). Then
\begin{eqnarray*}
g_1(\varp)-g_1(\psi)=\varp K(|\varp|^2)- \psi K(|\psi|^2)= (\varp-\psi)
K(|\varp|^2)+\psi K(|\varp|^2-|\psi|^2).
\end{eqnarray*}
Thus 
\begin{eqnarray*}
\|g(\varp)-g(\psi)\|_{L^{\rho'}}\leq\|(\varp-\psi)
K(|\varp|^2)\|_{L^{\rho'}} +\|\psi
K(|\varp|^2-|\psi|^2)\|_{L^{\rho'}}.
\end{eqnarray*}
Set $p=\frac{3r'\rho'}{3+2r'\rho'}$. As before from H\"older, Riesz and Sobolev inequalities we get 
\begin{eqnarray}
\|(\varp-\psi)K(|\varp|^2)\|_{L^{\rho'}}\leq
C\|\varp-\psi\|_{L^{r\rho'}}\|\varp\|_{\W}^{2}, 
\end{eqnarray}
and
\begin{eqnarray}
\|\psi
K(|\varp|^2-|\psi|^2)\|_{L^{\rho'}}&&\leq\|\psi\|_{L^{r\rho'}}\|K(|\varp|^2-|\psi|^2)\|_{L^{r'\rho'}}\nonumber\\
&&
\leq C\|\psi\|_{L^{r\rho'}}\|(\varp-\psi)(\varp+\psi)\|_{L^{p}}.
\end{eqnarray}

From the conditions imposed on $r,\rho$ we conclude that for
$\beta=\frac{r\rho'}{p}$ we have $\beta>1$ and $\beta' p\in
[2,6]$. Thus by the same arguments and observing that $\beta p=r\rho'$, we get
 
\begin{eqnarray}
\|(\varp-\psi)(\varp+\psi)\|_{L^{p}}&\leq&
\|(\varp-\psi)\|_{L^{\beta p}}\|(\varp+\psi)\|_{L^{\beta' p}}\leq
2M\soc\|\varp-\psi\|_{L^{r\rho'}}.
\end{eqnarray}

Finally putting all together, we get
\begin{eqnarray}
\|g_1(\varp)-g_1(\psi)\|_{L^{\rho'}}\leq \big( CM^{2}+2\soc
M^2\big)\|\varp-\psi\|_{L^{r\rho'}}.
\end{eqnarray}
 
and (ii) is proved.

\end{proof}

In order to show that the function $g_2$ satisfies similar conditions
as $g_1$, we shall need further auxiliary results  concerning the
function $g_1$.

\begin{lem} For every $\psi\in\W$, we have 
\begin{eqnarray} 
\|g_1(\psi)\|_{\Lb}^2\leq \rc\soc\soc'\|\psi\|_{\W}^4.
\end{eqnarray}
\label{p2}
\end{lem}

\begin{proof} It suffices to prove the inequality for positive functions from
  $\W$. Let $\psi$ be such a function and $\alp\geq 3$. Then $2<2\alp'\leq 6$,
  $1<\gamma:=\frac{6\alp}{3+4\alp}<3/2$ and $2\gamma\in[2,6]$. Thus
  by H\"older inequality together with Riesz and Sobolev inequalities we
  obtain
\begin{eqnarray} 
\|g_1(\psi)\|_{\Lb}^2&\leq&(\int\psi^{2\alp'})^{1/\alp'}(\int(K(\psi^2))^{2\alp})^{1/\alp}\leq
\soc\rc \|\psi\|_{\W}^2(\int\psi^{2\gamma}\,dx)^{1/\gamma}\nonumber\\
&&\leq\soc\soc'\rc\|\psi\|_{\W}^4.
\end{eqnarray}
\end{proof}
Here are the main properties of the function $g_2$
\begin{lem} The function, $g_2$ is well defined and satisfies:
\begin{itemize}
\item[i]) $g_2\in C\big(\W,\Wd)$.
\item[ii)] $g_2\in C\big(\W,L^2)$.
\item[iii)] For every $M>0$ and every $\varp,\psi\in\W$ such that
  $\|\varp\|_{\W}$, $\|\psi\|_{\W}\leq M$,  we have
$$ 
\|g_2(\varp)-g_2(\psi)\|_{L^{2}}\leq M^3A\|\varp-\psi\|_{L^{2}},
$$
where $A$ is a constant depending only on Riesz and Sobolev constants.
\end{itemize}
\label{lip-g2}
\end{lem}

\begin{proof} Clearly $g_2$ is well defined.\\
The property (i) follows from (ii), and the latter one follows from
(iii) if we just show that $g_2$ from $\W$ into $L^2$ is well
defined. On the other we observe that $g_2(\psi)=\psi G_1(\psi)$ and
the function $G_1:\W\to\R$ is well defined. Thus $g_2$ is well defined
as well from  $\W$ into $L^2$.\\
(iii) Let $M>0$ and  $\varp,\psi\in \W$ such that
$\|\varp\|_{\W}$, $\|\psi\|_{\W}\leq M$. Then 
\begin{eqnarray}
\|g_2(\varp)-g_2(\psi)\|_{L^{2}}\leq
G_1(\varp)\|\varp-\psi\|_{L^{2}}+|G_1(\varp)-G_1(\psi)|\|\psi\|_{\Lb}.
\end{eqnarray}
We now proceed to establish the sought estimate for each term of the letter
inequality separately.\\
For the first term of RHS, making use of Lemma\ref{p2}, we get  
\begin{eqnarray}
G_1(\varp)&=&\int|\varp|^2K(|\varp|^2)\,dx\leq\|\varp\|_{\Lb}\|g_1(\varp)\|_{\Lb}\leq
M^3(\rc\soc\soc')^{1/2}.
\end{eqnarray}

yielding

\begin{eqnarray}
G_1(\varp)\|\varp-\psi\|_{L^{2}}\leq M^3(\rc\soc\soc')^{1/2}\|\varp-\psi\|_{L^{2}}.
\label{part1}
\end{eqnarray}

To estimate the second term, we introduce the function 
\begin{eqnarray} 
\Gamma:=[0,1]\to\R,\ t\mapsto G_1\big(t\varp+(1-t)\psi\big).
\end{eqnarray}

Obviously $G_1$ is of class $C^1$ (Fr\'echet) and $G_1'=4g_1$. Thus $\Gamma$
is differentiable as well and
\begin{eqnarray}  
\Gamma'(t)=4{\rm Re}\int g_1(t\varp+(1-t)\psi)(\varp-\psi)\,dx.
\end{eqnarray}
Hence 
\begin{eqnarray}
|\Gamma(1)-\Gamma(0)|=&&|G_1(\varp)-G_1(\psi)|\leq 4\sup_{t\in
 [0,1]}\big|\int g_1(t\varp+(1-t)\psi)(\varp-\psi)\,dx\big|\nonumber\\
&& \leq 4\|\varp-\psi\|_{L^{2}} \sup_{t\in
 [0,1]}\big\|g_1(t\varp+(1-t)\psi)\big\|_{\Lb}\nonumber\\
&&
\leq 4(\rc\soc\soc')^{1/2}\|\varp-\psi\|_{L^{2}}(\|\varp\|_{\W}+\|\psi\|_{\W})^2\nonumber\\
&&
\leq 16M^2(\rc\soc\soc')^{1/2}\|\varp-\psi\|_{L^{2}}.
\label{part2}
\end{eqnarray}
 
Finally putting equations (\ref{part1}) and (\ref{part2}) together
yields the result, which completes the proof.

\end{proof}

Lastly we establish the most important feature of functions $g_1,g_2$, namely their  weak continuity.\\
In the sequel we denote by 
\begin{eqnarray}
\tilde{K}_n\phi:=\int_{\{|\cdot-y|\leq 1/n\}}\frac{|\phi(y)|^2}{|\cdot-y|}\,dy.
\end{eqnarray}
We shall also designate by $B$ various open balls in $\R^3$. 

\begin{lem} \begin{itemize}
\item[i)] The function $g_1:\W\to L^{\rho'}$ is continuous w.r.t. the weak
  topologies of $\W$ and $\Lbr$.
\item[ii)] The function $g_2:\W\to\Lb$ is continuous w.r.t. the weak topologies of
  $\W$ and $\Lb$.
\end{itemize}
\label{weak-cont}
\end{lem}

\begin{proof}

We will repeatedly use the known fact that for every open ball
$B\subset\R^3$ and every $2\leq s<6$ the space $\W$ embeds compactly into
$L^s(B)$.\\
(i): Let $(\psi_n)\subset\W$, $\psi\subset\W$ such that $\psi_n\rhup\psi$
in $\W$,  $B$ be an open ball in $\R^3$ and $w\in L^{\rho}$. Then
\begin{eqnarray}
\int
w\big(g_1(\psi_n)-g_1(\psi)\big)\,dx=\int(\psi_n-\psi)wK(|\psi_n^2|)\,dx+\int\psi
wK\big(|\psi_n|^2-|\psi|^2\big)\,dx
\end{eqnarray}
We decompose the first integral into
\begin{eqnarray*}
I_1(n):=\int(\psi_n-\psi)wK(|\psi_n^2|)\,dx=\int_B(\psi_n-\psi)wK(|\psi_n^2|)\,dx+\int_{B^c}(\psi_n-\psi)wK(|\psi_n^2|)\,dx
\end{eqnarray*}
Let $\rho$ and $r$ be the exponents given by Lemma\ref{lip-g1}. By
H\"older's inequality and the computations made in he proof of
Lemma\ref{lip-g1}, we get
\begin{eqnarray}
\int_B(\psi_n-\psi)wK(|\psi_n^2|)\,dx&&\leq\|w\|_{L^{\rho}}\big(\int_B|\psi_n-\psi|^{r\rho'}\,dx\big)^{1/{r{\rho'}}}\big(\int_B(|K(|\psi_n|^2)^{r'\rho'}\,dx\big)^{1/{r'\rho'}}\nonumber\\
&&
\leq C\|w\|_{L^{\rho}(B)}\|\psi_n\|_{\W}^2\big(\int_B|\psi_n-\psi|^{r\rho'}\,dx\big)^{1/{r{\rho'}}}\to 0,
\end{eqnarray}
by the fact that s $2\leq r\rho'<6$ and  the compactness of the embedding of $\W$ into $L^{r\rho'}$.\\ By the same arguments we get
\begin{eqnarray}
\int_{B^c}(\psi_n-\psi)wK(|\psi_n^2|)\,dx&&\leq\|w\|_{L^{\rho}(B^c)}\big(\int_{B^c}|\psi_n-\psi|^{r\rho'}\,dx\big)^{1/{r{\rho'}}}\big(\int_{B^c}(|K(|\psi_n|^2)^{r'\rho'}\,dx\big)^{1/{r'\rho'}}\nonumber\\
&&
\leq C\|w\|_{L^{\rho}(B^c)}\|\psi_n\|_{\W}^2\big(\int|\psi_n-\psi|^{r\rho'}\,dx\big)^{1/{r{\rho'}}}\nonumber\\
&&\leq C\|w\|_{L^{\rho}(B^c)}.
\end{eqnarray}
Now given $\epsilon>0$, we choose $B$ so that $\|w\|_{L^{\rho}(B^c)}<\epsilon$
and get
\begin{eqnarray}
 \int_{B^c}(\psi_n-\psi)wK(|\psi_n^2|)\,dx\leq \epsilon,\ \forall\,\epsilon>0,
\end{eqnarray}
yielding the convergence toward zero of $I_1(n)$.\\
We also decompose the second integral
\begin{eqnarray*}
I_2(n):=&&\int\psi wK\big(|\psi_n|^2-|\psi|^2\big)\,dx=
\int_BwK\big(|\psi_n|^2-|\psi|^2\big)\,dx +\int_{B^c}\psi
wK\big(|\psi_n|^2-|\psi|^2\big)\,dx\\
&&=\int\big(|\psi_n|^2-|\psi|^2\big)K(1_Bw\psi)\,dx+\int_{B^c}\psi
wK\big(|\psi_n|^2-|\psi|^2\big)\,dx,
\end{eqnarray*}
and rewrite the integral
\begin{eqnarray*}
I_2'(n):=\int\big(|\psi_n|^2-|\psi|^2\big)K(1_Bw\psi)\,dx&&=\int_B\big(|\psi_n|^2-|\psi|^2\big)K(1_Bw\psi)\,dx\\
&&
+\int_{B^c}\big(|\psi_n|^2-|\psi|^2\big)K(1_Bw\psi)\,dx.
\end{eqnarray*}
Choose $p$ so hat  $1\leq p<\frac{6\rho}{7\rho-6}$. Then
$1<p<\frac{3}{2}$. Setting $\alp=\frac{p}{p-1}$ and
$\beta=\frac{3\alp}{2\alp+3}$, yields 
\begin{eqnarray}
1<p<\frac{3}{2}\ {\rm and}\ 2\leq\frac{\rho\beta}{\rho-\beta}<6.
\label{p-beta}
\end{eqnarray} 
Using H\"older's inequality together with Riesz's potential estimate (\ref{riesz}) we obtain
\begin{eqnarray}
\int_B\big||\psi_n|^2-|\psi|^2\big|K(1_Bw\psi)\,dx&&\leq\||\psi_n|^2-|\psi|^2\|_{L^p(B)}\|K(w\psi)\|_{L^{\alp}(B)}\nonumber\\
&&
\leq C(M)\|\psi_n-\psi\|_{L^{2p}(B)}\|w\|_{L^{\rho}(B)}\|\psi\|_{L^{\frac{\rho\beta}{\rho-\beta}}(B)}.
\end{eqnarray}
Taking conditions (\ref{p-beta}) into account yields
$$  
\|\psi\|_{L^{\frac{\rho\beta}{\rho-\beta}}(B)}<\infty\ {\rm and}\ 
\lim_{n\to\infty}\|\psi_n-\psi\|_{L^{2p}(B)}=0.
$$
We conclude thereby that the latter integral tends to zero as $n\to\infty$.\\
By the same way we get the estimate
\begin{eqnarray}
\int_{B^c}\big||\psi_n|^2-|\psi|^2\big|K(1_Bw\psi)\,dx&&\leq
\||\psi_n|^2-|\psi|^2\|_{L^p(B^c)}\|K(1_Bw\psi)\|_{L^{\alp}(B^c)}\nonumber\\
&&\leq C(M)\|K(1_Bw\psi)\|_{L^{\alp}(B^c)}.
\end{eqnarray}
The already made calculus shows that 
\begin{eqnarray}
\int\big|K(1_Bw\psi)\big|^{\alp}\,dx<\infty.
\end{eqnarray}
Whence choosing, for every $\epsilon>0$, a $B$ such that
$\|K(1_Bw\psi)\|_{L^{\alp}(B^c)}<\epsilon$, gives
\begin{eqnarray}
\int_{B^c}\big||\psi_n|^2-|\psi|^2\big|K(1_Bw\psi)\,dx<\epsilon,\ \forall\ \epsilon>0,
\end{eqnarray}
and thereby $\lim_{n\to\infty}I_2'(n)=0$.
Finally writing
\begin{eqnarray*}
\int_{B^c}w\psi K\big(|\psi_n|^2-|\psi|^2\big)\,dx=\int\big(|\psi_n|^2-|\psi|^2\big)K(1_{B^c}w\psi )\,dx
\end{eqnarray*}
and using the same techniques as before we obtain 
\begin{eqnarray}
\lim_{n\to\infty}\int_{B^c}w\psi K\big(|\psi_n|^2-|\psi|^2\big)\,dx<\epsilon,\ \forall\ \epsilon>0.
\end{eqnarray}
Thus $\lim_{n\to\infty}(I_1(n)+I_2(n))=0$,,which was to be proved.\\
(ii): Let $(\psi_n),\psi$ be as before. We rewrite
\begin{eqnarray}  
g_2(\psi_n)&&=\psi_n(\psi_n,g_1(\psi_n))_{L^{\rho},\Lbr}=\psi_n\big(\psi_n,g_1(\psi_n)-g_1(\psi)\big)_{L^{\rho},\Lbr}\nonumber\\
&&
+\psi_n\big(\psi_n-\psi,g_1(\psi)\big)_{L^{\rho},\Lbr}+\psi_n\big(\psi,g_1(\psi)\big)_{L^{\rho},\Lbr}\nonumber\\
&&\rhup\psi(\psi,g_1(\psi))_{L^{\rho},\Lbr}=g_2(\psi),\ {\rm in}\ \Lb\ {\rm
  as}\ n\to\infty.
\end{eqnarray}

which finishes the proof.

\end{proof}

\section{Existence and uniqueness}

Thanks to the properties of $g_1,g_2$, we conclude that a function $\psi\in
L^{\infty}(I,\W)$ solves the (frNSE) if and only if it satisfies Duhamel's formula:
\begin{eqnarray}
\psi(t)={\rm e}^{i\alp_1tH}\varphi+i\alp_2\int_0^t{\rm
  e}^{i\alp_1(t-s)H}g_1(\psi(s))\,ds-i\alp_2\int_0^t{\rm
  e}^{i\alp_1(t-s)H}g_2(\psi(s))\,ds,\ \forall\,t\in I,
\label{duhamel}
\end{eqnarray}
where $H$ stands for Laplace operator on the Euclidean space
$\R^3$.\\
Observe that since $g_1$ does not preserve the space $\W$ it is not possible to use a
fixed point argument to solve the problem directly. Also the occurrence of
complex coefficients in the (frNSE) does not permit use of classical results,
especially those based on 'conservation laws' (see \cite{cazenave})to solve the problem. Instead we shall truncate
the Newton kernel, construct a sequence of approximate solutions then pass to
the limit.\\
To that end we introduce the sequence of functions $(f_n)$, $n\in\N^*$:
\begin{eqnarray*} 
f_n(\psi):=\psi\int_{\{|\cdot-y|>1/n\}}\frac{|\psi(y)|^2}{|\cdot-y|}\,dy:=\psi K_n(|\psi|^2).
\end{eqnarray*}

The function $f_n$ enjoys the following properties

\begin{lem}
\begin{itemize} 
\item[i)] For every $\psi\in\W$, $|f_n(\psi)|\leq |g_1(\psi)|$.
\item[ii)] For every $\psi\in\W$, $f_n(\psi)\in\W$.
\item[iii)] Let $\rho,r_1$ be the exponents given by Lemma\ref{lip-g1}. Then for
  every $0<M<\infty$, there is a constant $C(M)$ such that for each $n\in\N^*$,
$$
\|f_n(\varp)-f_n(\psi)\|_{L^{\rho'}}\leq C(M)\|\varp-\psi\|_{L^{r_1}},
$$  
for every $\varp,\psi\in\W$ such that $\|\varp\|_{\W}\leq M,\
\|\psi\|_{\W}\leq M$.
\end{itemize} 
\label{fn}
\end{lem} 

\begin{proof} The proof of property (i) is obvious. To prove (ii) observe that
  by (i) since for every $\psi\in L^2,\ g_1(\psi)\in L^2$ then $f_n(\psi)\in L^2$ as
  well.  Now a direct computation yields  that for every $\psi\in\W$, $\nabla
  f_n(\psi)=K_n(|\psi|^2)\nabla\psi+\psi\nabla K_n(|\psi|^2)$ and
\begin{eqnarray}
|K_n(|\psi|^2)\nabla\psi|\leq n\|\psi\|_{L^2}^2|\nabla\psi|\in L^2,
\end{eqnarray}
\begin{eqnarray}
|\psi||\nabla
 K_n(|\psi|^2)|=|\psi|\big|\int_{\{|\cdot-y|>1/n\}}\frac{(\cdot-y)|\psi(y)|^2}{|\cdot-y|^3}\,dy\big|\leq
 n^2|\psi|\|\psi\|_{L^2}^2\in L^2.
\end{eqnarray}
The proof of (iii) follows by using (i), the fat that $\big|K_n\phi\big|\leq
K(|\phi|)$ and Lemma\ref{lip-g1}.

\end{proof}

The most important statement of the latter lemma is property (iii), which
indicates that the Lipschitz constant, as well as the exponents $\rho$ and $r_1$, are independent of the integer $n$.\\  
Consider now the approximate problem
\begin{eqnarray}
\begin{cases}
i\frac{\partial\psi}{\partial t}=-\alp_1\Delta\psi+i\alp_2f_n(\psi)-i\alp_2g_2(\psi) & \\
\psi(0)=\varphi &  \end{cases}.
\label{app-frnse}
\end{eqnarray} 

As before, we assert that $\psi_n\in L^{\infty}(I,\W)$ solves
(\ref{app-frnse}) if and only if it satisfies the related Duhamel's formula

\begin{eqnarray}
\psi_n(t)={\rm e}^{i\alp_1tH}\varphi+i\alp_2\int_0^t{\rm
  e}^{i\alp_1(t-s)H}f_n(\psi_n(s))\,ds-i\alp_2\int_0^t{\rm
  e}^{i\alp_1(t-s)H}g_2(\psi_n(s))\,ds.
\label{duhameln}
\end{eqnarray}
 
Thanks to the already observed fact that both $f_n$ and $g_2$ preserve the space
$\W$ it is possible to solve the latter equation via a fixed
point-argument. However, we shall use that argument  not for
the operator defined by the RHS of (\ref{duhameln}), but for some power
of it. To this end we  use a theorem due to Weissinger (see \cite{heuser}):

\begin{theo} (Weissinger) Let $Y$ be a Banach space and
  $(\alp_k)$ a sequence of positive numbers  such that $\sum\alp_k<\infty$. Let
  $F\subset Y$ be closed and $A:F\to F$ be an operator such that
$$
\|A^k\varp-A^k\psi\|\leq\alp_k\|\varp-\psi\| ,\ \forall\,k\in\N\ {\rm and}\,\forall\,
\varp,\psi\in Y.
$$
Then $A$ possesses a unique fixed point. Furthermore the fixed point can be
obtained as the limit of the sequence defined by $\psi_0=\Phi\in F$,
and $\psi_{k+1}=A\psi_k$. 
\end{theo}
We shall also make use of the known fact that the operator ${\rm
  e}^{itH}$ is unitary on each of the spaces $\Lb$, $\W$ and $\Wd$.\\

Now let $M>0$, $T>0$ and $\varp\in\W$ be given. Set $I=[-T,T]$, $X:=L^{\infty}(I,\W)$ and 
$$
F_M:=\{\psi\in X:\,\|\psi-{\rm e}^{i\cdot H}\varp\|_X\leq M\}.
$$

We will first determine $T$ so that for every integer $n$ the operators
\begin{eqnarray*} 
A_n:F_M\to X,\,\psi\mapsto{\rm e}^{i\alp_1tH}\varphi+i\alp_2\int_0^t{\rm
  e}^{i\alp_1(t-s)H}f_n(\psi(s))\,ds-i\alp_2\int_0^t{\rm
  e}^{i\alp_1(t-s)H}g_2(\psi(s))\,ds
\end{eqnarray*}
maps $F_M$ into itself.\\
For $\psi\in X$ set 
\begin{eqnarray}
S_n\psi:=\int_0^t{\rm e}^{i\alp_1(t-s)H}f_n(\psi(s))\,ds\ {\rm and}\ U\psi:=\int_0^t{\rm
  e}^{i\alp_1(t-s)H}g_2(\psi(s))\,ds.
\end{eqnarray}

We still designate par $\rho$ and $r_1$ the exponents  given by Lemma\ref{lip-g1}
and we will omit in the notation  the  dependence of $F$ on $M$. Various
constants depending on $M$ will be denoted by $C(M)$.

\begin{lem} Let $M>0$ and $\varp\in\W$ such that $\|\varp\|_{\W}\leq M$. There
  is $0<T=T(M)$ such that for every $n\in\N^*$ the operator  $A_n$ maps $F$ into itself.
\label{contraction}
\end{lem}

\begin{proof} Let $M>0$ be fixed and $\psi\in F_M$. Then

\begin{eqnarray*}
\|S_n\psi\|_{\W}+\|U\psi\|_{\W}&=&\|S_n\psi\|_{\Wd}+\|U\psi\|_{\W}\\
&&
\leq\int_0^{|t|}\|{\rm e}^{i\alp_1(t-s)H}\|_{\Wd,\Wd}\|f_n(\psi(s))\|_{\Wd}d\,s\\
&& 
+ \int_0^{|t|}\|{\rm e}^{i\alp_1(t-s)H}\|_{\W,\W}\|g_2(\psi(s))\|_{\W}d\,s\\
&&
\leq C_1\int_0^{|t|}\|f_n(\psi(s))\|_{L^{\rho'}}d\,s+C_2\int_0^{|t|}\|\psi(s)\|_{\W}^{4}d\,s\\
&&
\leq
C_1(M)\int_0^{|t|}\|g_1(\psi(s))\|_{L^{\rho'}}d\,s+C_2(M)\int_0^{|t|}\|\psi(s)\|_{\W}^{4}d\,s\\
&&
\leq C_1(M)\int_0^{|t|}\|\psi(s)\|_{L^{r_1}}d\,s+C_2(M)\int_0^{|t|}\|\psi(s)\|_{\W}^{4}d\,s\\
&&
\leq C(M)|t|\big(\|\psi\|_{X}+\|\psi\|_{X}^{4}\big),
\end{eqnarray*}
yielding, for $0<|t|\leq T$
\begin{eqnarray}
\|A_n\psi-{\rm e}^{i\cdot H}\varp\|_X& &\leq
\alp_2\|S_n\psi\|_{X}+\alp_2\|U\psi\|_{X}\\
&&\leq C(M)T(M+\|\varp\|_{\W})(1+ (M+\|\varp\|_{\W})^3).
\end{eqnarray}
Finally we choose $T$ small  so that 
$$
 2MC(M)(1+ 8M^3)T\leq M,
$$
which completes the proof.
\end{proof}

Now we proceed to show that for each integer $n$, operators $A_n$ satisfy the conditions
demanded by Weissinger's theorem.

\begin{lem} Let $M>0, T>0$ and $\varp\in\W$, $\|\varp\|_{\W}\leq M$ be fixed. Then there is a
  constant $C$ depending only on $M$, Riesz's and Sobolev constants such that
  for every $t\in[-T,T]$, every integer $k,n$ and every $\psi_1,\psi_2\in
  X$ such that $\|\psi_1\|_X,\|\psi_2\|_X\leq M$   we have
\begin{eqnarray}   
\|A_n^k\psi_1(t)-A_n^k\psi_2(t)\|_{\W}\leq\frac{(Ct)^k}{k!}\|\psi_1-\psi_2\|_X.
\end{eqnarray}
\label{exponentiel}
\end{lem}

\begin{proof} The proof runs by induction, with the help of the local
  Lipschitz property of both functions $g_1,g_2$.\\
We will only give the idea how to get the estimate for the $k=1$. For general
  $k$ the estimate follows by direct induction.\\
Without loss of generality we assume that $t\geq 0$. Let $\psi_1,\psi_2\in
  X$ be such that $\|\psi_1\|_X\leq M$, $\|\psi_2\|_X\leq M$. Let $r_1$ be as given by
  Lemma\ref{lip-g1}. Then 

 \begin{eqnarray*}   
\|S_n\psi_1(t)-S_n\psi_2(t)\|_{\W}&=&\|S_n\psi_1(t)-S_n\psi_2(t)\|_{\Wd}\\
&&\leq\int_0^t\|{\rm
  e}^{i\alp_1(t-s)H}\|_{\Wd,\Wd}\|f_n(\psi_1(s))-f_n(\psi_2(s))\|_{\Wd}d\,s\\
&&
 \leq \int_0^t \|f_n(\psi_1(s))-f_n(\psi_2(s))\|_{L^{\rho'}}d\,s\\
&&\leq C\int_0^t 
  \|g_1(\psi_1(s))-g_1(\psi_2(s))\|_{L^{\rho'}}d\,s\nonumber\\
&&
\leq C(M)\int_0^t\|g_1(\psi_1(s))-g_1(\psi_2(s))\|_{L^{r_1}}d\,s\nonumber\\
&&
\leq C(M)\int_0^t \|\psi_1(s)-\psi_2(s)\|_{\W}d\,s\\
&&
\leq tC(M)\|\psi_1-\psi_2\|_X.
\end{eqnarray*}

Here we used the fact that$L^{\rho'}$ embeds continuously into $\Wd$  . By the same ideas we achieve

\begin{eqnarray*}   
&\|U\psi_1(t)-U\psi_2(t)\|_{\W}&\leq\int_0^t\|{\rm
  e}^{i\alp_1(t-s)H}\|_{\W,\W}\|g_2(\psi_1(s))-g_2(\psi_2)\|_{\W}d\,s\nonumber\\
&&
\leq C(M)\int_0^t \|\psi_1(s)-\psi_2(s)\|_{\W}d\,s\leq C(M)t\|\psi_1-\psi_2\|_X.
\end{eqnarray*}

Thus

\begin{eqnarray}   
\|A_n\psi_1(t)-A_n\psi_2(t)\|_{\W}\leq C(M)t\|\psi_1-\psi_2\|_X.
\label{A2}
\end{eqnarray}

\end{proof}

We are now in position to affirm the local solvability of the approximate
problem (\ref{app-frnse}).

\begin{theo} Let $M>0$ and $\varphi\in\W$ such that $\|\varp\|_{\W}\leq M$ be
  fixed. Then there is $T_M>0$ such that for every $n\in\N^*$ problem (\ref{app-frnse}) has a unique solution, $\psi_n$, in the space $L^{\infty}([-T_M,T_M],\W)$.
Further the solution may be gained as the limit of the sequence $\psi_0=\Phi$,
$\psi_{k+1}=A_n\psi_k$, where $\Phi$ is any element from $L^{\infty}([-T_M,T_M],\W)$.  
\label{app-existence}   
\end{theo}

\begin{proof} Making use of Duhamel's formula, we have simply to check that assumptions of Weissinger's theorem
  are fulfilled.\\
Let $M>0$ and $\varp\in\W$, $\|\varp\|_{\W}\leq M$ be fixed. By
Lemma\ref{contraction}, for every $n\in\N^*$, there is $T:=T_M>0$ such that operators $A_n$ map the closed ball of $X:=L^{\infty}([-T,T],\W)$ of radius $M$
  and centered on ${\rm e}^{itH}\varp$, $F$ into itself.\\
Setting $\beta_k:=\frac{(C(M)T)^k}{k!}$, $\forall\,k\in\N$, we obtain by  Lemma\ref{exponentiel}
\begin{eqnarray}   
\|A_n^k\psi_1-A_n^k\psi_2\|_{X}\leq \beta_k\|\psi_1-\psi_2\|_X,\ \forall\,k\in\N,\,n\in\N^*
\end{eqnarray}
with $\sum_{k=0}^{\infty}\beta_k={\rm exp}(C(M)T)$, which completes the proof.
\end{proof}

For our later purposes we establish continuous dependence of the approximate
solution (solution of the approximate problem) w.r.t. the initial data. 

\begin{prop}
Let $0<M'<M$ and
  $\varp,\tilde\varp\in\W$ be such that $\|\varp-\tilde\varp\|_{\W}\leq
  M-M'$. Set $T:=\min(T_M(\varp),T_M(\tilde\varp))$ and $\psi_n$,
  resp. $\tilde\psi_n$ the local solution of (\ref{app-frnse}) with
  $\psi_n(0)=\varp$, resp. $\tilde\psi_n(0)=\tilde\varp$. Then there is a
  constant $C$ such that
\begin{eqnarray}
\sup_{|t|\leq T}\|\psi_n(t)-\tilde\psi_n(t)\|_{\W}\leq
C\|\varp-\tilde\varp\|_{\W},\ \forall\,n\in\N^*.
\end{eqnarray}
\label{app-continuity}
\end{prop}

\begin{proof} As observed in Theorem\ref {existence}, the solution
  $\psi_n$ on $[-T,T]$ is given as the limit of the sequence
\begin{eqnarray} 
\Phi_0\in\{u\in L^{\infty}([-T,T],\W):\,\|u-e^{i\cdot
  H}\varp\|_{L^{\infty}([-T,T],\W)}\leq M\},\ \Phi_{k+1}=A_n\Phi_k.
\end{eqnarray}
By the conditions imposed on $M,M',\varp$ and $\tilde\varp$ we have
 $\|\tilde\psi_n-e^{i\cdot
  H}\varp\|_{L^{\infty}([-T,T],\W)}\leq M$. Thus we can choose
  $\Phi_0=\tilde\psi_n$ on $[-T,T]$.\\
Setting $h_n:=i\alp_2(f_n-g_2)$ we get: $\forall\,t\in [-T,T]$,
\begin{eqnarray*}
\tilde\psi_n(t)=e^{i\alp_1tH}\tilde\varp+\int_0^te^{\alp_1i(t-s)H}h_n(\tilde\psi_n(s))\,ds=e^{i\alp_1 tH}\tilde\varp+\int_0^te^{\alp_1i(t-s)H}h_n(\Phi_0(s))\,ds.
\end{eqnarray*}
Thus, for every $t\in [-T,T]$, we have
\begin{eqnarray}
\Phi_1(t)-\Phi_0(t)=\Phi_1(t)-\tilde\psi_n(t)=e^{i\alp_1
  tH}(\varp-\tilde\varp),
\end{eqnarray}
and 
\begin{eqnarray}
\sup_{|t|\leq
  T}\|\Phi_1(t)-\Phi_0(t)\|_{\W}\leq\|\varp-\tilde\varp\|_{\W}.
\end{eqnarray}
On the other hand we have, for every $n,k\in\N$
\begin{eqnarray}
\psi_n-\Phi_k=\sum_{j=k+1}^{\infty}(\Phi_{k+1}-\Phi_k),
\end{eqnarray}
yielding for $n=1$and for $\beta_k,C$  as given in the
 proof of Theorem\ref{existence}
\begin{eqnarray}
\|\psi_n-\Phi_1\|_{L^{\infty}([-T,T],\W)}&&\leq\big(\sum_{k=2}^{\infty}\beta_k\big)\|\Phi_1-\Phi_0\|_{L^{\infty}([-T,T],\W)}\nonumber\\
&&\leq{\rm exp}(CT)\|\varp-\tilde\varp\|_{\W}.
\end{eqnarray}
Putting all together gives
\begin{eqnarray}
\|\psi_n-\tilde\psi_n\|_{L^{\infty}([-T,T],\W)}\leq{\rm
  exp}(CT)\|\varp-\tilde\varp\|_{\W},\ \forall,\ n\in\N
\end{eqnarray}
and the proof is finished.

\end{proof}

We stress that the constant occurring in the estimate given by
Proposition\ref{continuity} does not depend on $n$ but only on $M$ and $M'$.\\

Next we shall rely on the Theorem\ref{app-existence}  result to prove local existence of
solutions for (\ref{frnse}).

\begin{theo} Let $M>0$ and $\varphi\in\W,\ \|\varphi\|_{\W}\leq M$ be fixed. Then there is
  $T_M>0$ such that problem (\ref{frnse}) has a solution, $\psi\in L^{\infty}([-T_M,T_M],\W)$.
\label{existence}   
\end{theo}

\begin{proof}  On the light of Theorem\ref{app-existence}, there is $T:=T_M$
  and a sequence of approximate solutions $(\psi_n)\subset
  X:=L^{\infty}\big([-T,T],\W\big)$ of problem (\ref{app-frnse}). Thus, for
  every $n$, $\psi_n$ satisfies

\begin{eqnarray}
i\frac{\partial\psi_n}{\partial
  t}=-\alp_1\Delta\psi_n+i\alp_2f_n(\psi_n)-i\alp_2g_2(\psi_n),\ {\rm in}\ \Wd.
\label{dist-eq}
\end{eqnarray} 
Making use of the uniform boundedness of
$(\psi_n)$ in $X$, we achieve
\begin{eqnarray}
\|\frac{\partial\psi_n}{\partial t}\|_{\Wd}&&\leq
C\big(\|\nabla\psi_n\|_{\W}+\|f_n(\psi_n)\|_{\Wd}+\|g_2(\psi_n)\|_{\Wd}\big)\nonumber\\
&&
\leq C(M)\big(1+\|f_n(\psi_n)\|_{\Lbr}+\|g_2(\psi_n)\|_{\W}\big)\nonumber\\
&&
\leq C(M)\big(1+\|f_n(\psi_n)\|_{\Lbr}+1)\leq
C(M)\big(2+\|g_1(\psi_n)\|_{\Lbr}\big)\nonumber\\
&&\leq C(M).
\end{eqnarray}
Therefore the sequence $(\psi_n)$ is uniformly bounded in
$$
Y:=X\cap W^{1,\infty}\big([-T,T],\Wd\big). 
$$
Thus (see \cite[ Proposition.1.3.14]{cazenave}) there is $\psi\in Y$ and a subsequence which
we denote also by  $(\psi_n)$ such hat
\begin{eqnarray}
\psi_n(t)\rightharpoonup\psi(t)\ {\rm in}\ \W,\,\forall t\in [-T,T].
\end{eqnarray}
Thus $\psi(0)=\varp$.\\
Let $\tilde K_n$, be the operators defined by
\begin{eqnarray}
\tilde K_n\phi:=\int_{|\cdot-y|<1/n}\frac{\phi(y)}{|\cdot-y|}\,dy,\ n\in\N^*.
\end{eqnarray}
Having Duhamel's formula (for $\psi_n$'s) in hand and rewriting
\begin{eqnarray}
f_n(\psi_n)=\psi_nK_n(|\psi_n|^2)=\psi_n(K-\tilde{K}_n)(|\psi_n^2|)
=g_1(\psi_n)-\psi_n\tilde{K}_n(|\psi_n^2|),
\end{eqnarray}
we get that for every $\phi\in\Lb\cap\Lbr$,
\begin{eqnarray}
<\psi_n(t),\phi>_{L^{\rho},\Lbr}&&=<{\rm e}^{i\alp_1tH}\varp,\phi>_{L^{\rho},\Lbr}+i\alp_2\int_0^t<{\rm
  e}^{i\alp_1(t-s)H}g_1(\psi_n(s)),\phi>_{L^{\rho},\Lbr}\,ds\nonumber\\
&&
-i\alp_2\int_0^t<{\rm
  e}^{i\alp_1(t-s)H}\big(\psi_n(s)\tilde
  K_n(|\psi_n^2|)(s)\big),\phi>_{L^{\rho},\Lbr}\,ds\nonumber\\
&&
-i\alp_2\int_0^t<{\rm
  e}^{i\alp_1(t-s)H}g_2(\psi_n(s)),\phi>_{L^{\rho},\Lbr}\,ds.
\end{eqnarray}
We claim that that $\|\tilde{K}_n\|_{L^p,L^p}\to 0$ for every
$1<p<\infty$. Indeed: For every $\phi\in L^p$, setting $q$ the conjugate
exponent of $p$ and
$$
G_n(x,y):=1_{\{|x-y|<1/n\}}|x-y|^{-1},
$$
we get
\begin{eqnarray} 
|\tilde K_n\phi(x)|\leq
 \big(C_n:=\sup_x\int G_n(x,y)\,dy\big)^{1/q}\int G_n(x,y)|\phi(y)|^p\,dy
\end{eqnarray}
and thereby
\begin{eqnarray}
\|\tilde K_n\|_{L^p,L^p}\leq C_n={c}/{n^2}\to 0\ {\rm as}\ n\to\infty.
\end{eqnarray}
Thus we get by  Proposition\ref{weak-cont} and  use of the fact that ${\rm e}^{itH}$ maps continuously $\Lbr$ into
$L^{\rho}$ for every $t\neq 0$, together with dominated convergence theorem, that for every
$\phi\in\Lb\cap\Lbr$,
\begin{eqnarray}
<\psi_n(t),\phi>_{L^{\rho},\Lbr}\to&&<{\rm e}^{i\alp_1tH}\varp,\phi>_{L^{\rho},\Lbr}+i\alp_2\int_0^t<{\rm
  e}^{i\alp_1(t-s)H}g_1(\psi(s)),\phi>_{L^{\rho},\Lbr}\,ds\nonumber\\
&&
-i\alp_2\int_0^t<{\rm
  e}^{i\alp_1(t-s)H}g_2(\psi(s)),\phi>_{L^{\rho},\Lbr}\,ds\nonumber\\
&&
=<\psi(t),\phi>_{L^{\rho},\Lbr},
\end{eqnarray}
yielding therefore
\begin{eqnarray}
\psi(t)={\rm e}^{i\alp_1tH}\phi+i\alp_2\int_0^t{\rm
  e}^{i\alp_1(t-s)H}g_1(\psi(s))\,ds-i\alp_2\int_0^t{\rm
  e}^{i\alp_1(t-s)H}g_2(\psi(s))\,ds.
\end{eqnarray}
Whence $\psi$ satisfies
\begin{eqnarray}
i\frac{\partial\psi}{\partial
  t}=-\alp_1\Delta\psi+i\alp_2f-i\alp_2g,\ {\rm in}\ \Wd,
\end{eqnarray} 
and $\psi$ is a solution of equation (\ref{frnse}).\\
{\em Uniqueness}: Follows from \cite[ Proposition4.2.3, p.85.]{cazenave}.

\end{proof}

\begin{prop} ({\em Blow-up alternative}) The blow-up alternative holds true for
  the solution of the (frNSE).
\label{blow-up}
\end{prop} 
The proof is quite standard so we omit it.

Yet we will describe how does the local solution behaves w.r.t. the
initial data.

\begin{prop}({\em continuous dependence}) Let $(\varp_k)_k\subset\W$ and
  $\varp\subset\W$ be such that $\|\varp_k-\varp\|_{\W}\to 0$. Set $\tilde\psi_{k}$ resp. $\psi$ the
  local solution of the frictional Newton-Schr\"odinger equation with initial
  data $\varp_k$, resp. $\varp$. Then  there is $T>0$ such that $\lim_{k\to\infty}\|\psi-\tilde\psi_k\|_{L^{\infty}((-T,T),\W)}=0$.
\label{continuity}
\end{prop}

\begin{proof} Set $\psi_{n,k}$, resp. $\psi_n$ the solution of the approximate
  Newton-Schr\"odinger equation with initial data $\varp_k$,
  resp. $\varp$. Making use of Proposition \ref{app-continuity}, there are
 constants  $C,T>0$ depending only on $\|\varp\|_{\W}$ such that for large $k$
\begin{eqnarray}
\sup_{t\in[-T,T]}\|\psi_n(t)-\psi_{n,k}(t)\|_{\W}\leq C \|\varp_k-\varp\|_{\W},\ \forall\,n\in\N.
\end{eqnarray}
By the proof of the existence of Theorem\ref {existence} together with the
  uniqueness we conclude that for large $k$
\begin{eqnarray}
\psi_n\rhup\psi,\ \psi_{n,k}\rhup\tilde\psi_k,\ \in\W,\ \forall\,t\in[-T,T].
\end{eqnarray}
Whence by the weak lower semi continuity of the norm we get for large $k$
\begin{eqnarray}
\|\psi_n(t)-\tilde\psi_{k}(t)\|_{\W}&&\leq\liminf_{n\to\infty}\|\psi_n(t)-\psi_{n,k}(t)\|_{\W}\nonumber\\
&&
\leq C\|\varp_k-\varp\|_{\W}\ \forall\,t\in[-T,T],
\end{eqnarray}
yielding the result.

\end{proof}

\section{Concluding remarks}

We would like to stress that our method (except maybe for the proof
of uniqueness) still works in a general domain of $\R^3$. However, if
$\Omega\subset\R^3$ is bounded then, thanks to the properties of the Newton kernel
on bounded subsets, it is possible to use an $L^2$-Gronwall-type inequality to get
the uniqueness.\\
At this stage, we mention that our method suggets an abstract framework for
solving evolution equations related to some classes of positive operators.\\
Finally, we indicate some open problems related to the
(frNSE). The first one is, of course, that dealing with the global existence
of the solution. Here we expect that a global solution would exists provided the
energy of the initial data is small enough. We are yet working in this
direction. Furthermore if a global solution exits it is interesting to ask
about its large time behavior. For the (NSE), this question was already
investigated by Wada \cite{takeshi01}.\\
The second one is much more complicated: Having the frictional
Newton-Schr\"odinger equation proposed by Diosi \cite{diosi07} (which is still
unsolved to our best knowledge!)  in mind, one is tempted to replace Diosi's
kernel by an other one, say
\begin{eqnarray}     
N(x,y)=\int\int G(x,z)G(y,z')\,d\mu(z)\,d\mu(z'),
\end{eqnarray}
where $G$ is positive, symmetric and  $\mu$ is a positive Radon measure. The
immediate question that arises is under which conditions on the measure $\mu$
and on $G$ has the related (frNSE)  a solution(s)? Is it local or global
and is the related (frNSE) well-posed?  
In this stage, to illuminate the way, one has first to look for the problem with the kernel
proposed by Diosi.\\
The last problem is the obvious generalization of the above questions in
higher dimensions.\\

{\sc Acknowledgment}. The first author would like to thank the 'University of
Bielefeld', where parts of this joint work was elaborated, for the warm hospitality.

\bibliography{biblio-frnse}

\end{document}